\documentclass{article}

\usepackage{amsmath,amssymb,amsthm}
\usepackage{tikz}
\usepackage{color}
\usepackage[toc]{appendix}
\usepackage{graphicx}
\usepackage{fancyhdr}
\usepackage{enumitem}
\usepackage{bbm}
\usepackage{parskip}
\usepackage{float}
\usepackage{chngpage}
\usepackage{calc}
\usepackage{bigints}
\usepackage{array}
\usepackage{booktabs}
\usepackage{rotating}
\usepackage{multirow}
\usepackage{adjustbox}
\usepackage{tabularx}
\usepackage{verbatim}
\usepackage{mathtools}
\usepackage{ragged2e}
\usepackage[makeroom]{cancel}
\usepackage{caption}
\usepackage{hyperref}
\usepackage{caption}
\usepackage{subcaption}
\usepackage{appendix}
\usepackage{pgfplots}
\pgfplotsset{compat=1.16}
\usepackage[english]{babel}
\usepackage{hyphenat}
\usepackage[makeindex]{imakeidx}
\usepackage{xcolor}

\usepackage{amsmath}

\usetikzlibrary{datavisualization}
\usetikzlibrary{matrix}
\usetikzlibrary{datavisualization.formats.functions}

\newtheorem{theorem}{Theorem}[section]

\newtheorem{corollary}[theorem]{Corollary}

\newtheorem{definition}[theorem]{Definition}

\newtheorem{remark}[theorem]{Remark}

\makeatletter
\def\section{\@startsection {section}{1}{\z@}{3.25ex plus 1ex minus
		.2ex}{1.5ex plus .2ex}{\large\bf}}
\def\subsection{\@startsection{subsection}{2}{\z@}{3.25ex plus 1ex minus
		.2ex}{1.5ex plus .2ex}{\normalsize\bf}}
\@addtoreset{equation}{section} 
\makeatother

\title{A note on Tricomi-type partial differential equations with white noise initial condition}
\author{Enrico Bernardi\thanks{Dipartimento di Scienze Statistiche Paolo Fortunati, Università di Bologna, Bologna, Italy. \\ \textbf{e-mail}: enrico.bernardi@unibo.it \textbf{ORCID}: 0000-0003-3923-1407} \and Alberto Lanconelli\thanks{Dipartimento di Scienze Statistiche Paolo Fortunati, Università di Bologna, Bologna, Italy.\\ \textbf{e-mail}: alberto.lanconelli2@unibo.it  \textbf{ORCID}: 0000-0001-8248-2151}}
\date{\today}

\begin{document}
	
	\maketitle
	
	\bigskip
	
	\begin{abstract}
		We study a class of Tricomi-type partial differential
		equations previously investigated in \cite{Y}. Firstly, we
		generalize the representation formula for the solution
		obtained there by allowing the coefficient in front of the second-order partial derivative
		with respect to $x$ to be any
		non integer power of $t$. Then, we analyze the robustness of
		that solution by taking the initial data to be  Gaussian
		white noise and we discover that the existence of a well-defined random field solution is lost upon the introduction of
		lower-order terms in the operator.  This phenomenon shows that,
		even though the Tricomi-type operators with or without lower-order terms are the same from the point of view of the theory
		of hyperbolic operators with double characteristics, their
		corresponding random versions exhibit different well posedness
		properties.  We also prove that for more regular initial data \textemdash
		specifically fractional Gaussian white noise with Hurst
		parameter $H\in (1/2,1)$\textemdash  the well posedness of the Cauchy problem for the Tricomi-type operator with lower-order term is restored. 
	\end{abstract}
	
	
	
	\allowdisplaybreaks
	
	\bigskip

\section{Introduction and statement of the main results}

The Tricomi equation \cite{tricomi} is a second-order partial
differential equation of mixed elliptic-hyperbolic type with the form:
\begin{align}\label{tricomi intro}
	u_{tt}(t,x)=tu_{xx}(t,x),\quad t,x\in\mathbb{R}.
\end{align}
The equation is hyperbolic in the half plane $t > 0$, elliptic in the half plane $t < 0$, and degenerates on the line $t = 0$. Many important problems in fluid mechanics and differential geometry can be reduced to corresponding problems for the Tricomi equation, particularly \emph{transonic flow} problems \cite{manwell} and \emph{isometric embedding} problems \cite{Qing} (see also \cite{tricomiintro} and the references quoted there). The Tricomi equation is a prototype of the Tricomi-type equation:
\begin{align}\label{gen intro}
	u_{tt}(t,x)=f(t)u_{xx}(t,x),\quad t,x\in\mathbb{R},
\end{align}
where $f:\mathbb{R}\to\mathbb{R}$ is a continuous function. 

In the series of papers \cite{Gelfand1,Gelfand2,Gelfand2bis,Gelfand3} the authors have investigated the fundamental solution of the Tricomi operator \eqref{tricomi intro}. In \cite{Barros-Neto} the method of partial Fourier transform is used to find explicit formulas for fundamental solutions for the following Tricomi-type operator in arbitrary space dimensions:
\begin{align}\label{tricomi intro multi}
	u_{tt}(t,x)=t\Delta_x u(t,x),\quad (t,x)\in\mathbb{R}^{n+1}.
\end{align}
More recently, the findings of \cite{Barros-Neto} have been generalized in \cite{Y} to the case of operators of the form
\begin{align}\label{tricomi intro multi 2}
	u_{tt}(t,x)=t^{2k}\Delta_x u(t,x),\quad (t,x)\in\mathbb{R}^{n+1},k\in\mathbb{N};
\end{align}
the author also provides representation formulas for the solution to related Cauchy problems. See also \cite{BealsKannai1,BealsKannai2} for explicit representations for the propagator of some generalized Tricomi-type operators in one space dimension.

The aim of the present note is twofold. We first propose a
probabilistic representation for the solution to the Cauchy
problem for the Tricomi-type equation \eqref{gen intro} with
$f(t)=t^{\alpha}$, $\alpha\geq 0$, thus yielding a
generalization (in one space dimension) of the results of
\cite{Y} for the operator \eqref{tricomi intro multi 2}. Then,
we exploit such representation to study the robustness of
solutions with respect to random initial conditions and the
introduction of lower-order terms of the type $u_x$. More
specifically, we set the initial data to be a Gaussian white
noise and investigate the existence of a well defined random
field solution for that problem. Notice that Cauchy problems
with Gaussian white noise initial conditions can be considered
by virtue of Duhamel's principle \cite{Evans} as an
intermediate step in the construction of mild solutions for
stochastic partial differential equations with additive
space-time Gaussian white noise, see
e.g. \cite{DaPrato,Rocknerbook,Pardouxbook}. We discover that
our Tricomi-type operator without lower-order term admits a
well defined Gaussian random field solution, while this does
not happen when a lower-order term $u_x$ is added. We also
show that such interference of the lower-order terms in the
well posedness of the Cauchy problem is absent when the
Tricomi-type operator reduces to the classical wave
operator. We complete our study pointing out that for smoother
initial data, i.e. fractional  Gaussian white noise with Hurst
parameter $H\in (1/2,1)$, the Cauchy problem for the
Tricomi-type operator with lower-order term  becomes well posed. 

We begin with our first main result whose proof is postponed
to Section \ref{proof 1}. Notice that the class of equations
considered here includes the family of problems investigated
in \cite{Y} ($\alpha=2k$, $k\geq 1$). From this point of view
Theorem \ref{th representation general} generalizes  Theorem 3.1 and Corollary 3.3 from  \cite{Y}.

\begin{theorem}\label{th representation general}
	Let $\alpha\geq 0$ and
	$\varphi:\mathbb{R}\to\mathbb{R}$ be twice
	continuously differentiable.  Then, the unique classical solution of the Cauchy problem for the Tricomi-type partial differential equation  
	\begin{align}
		\begin{cases}\label{tricomi}
			u_{tt}(t,x)=t^{\alpha}u_{xx}(t,x),& t>0,x\in\mathbb{R};\\
			u(0,x)=0, u_t(0,x)=\varphi(x),& x\in\mathbb{R},
		\end{cases}
	\end{align}
	can be represented as
	\begin{align}\label{representation}
		u(t,x)=t\mathbb{E}[\varphi(x+\xi(t) Z Y)],\quad t\geq 0,x\in\mathbb{R},
	\end{align}
	where:
	\begin{itemize}
		\item $Z$ is a Rademacher random variable, i.e. $\mathbb{P}(Z=1)=\mathbb{P}(Z=-1)=\frac{1}{2}$;
		\item $Y$ is a continuous random variable with probability density function
		\begin{align}
			p_Y(y)=
			\begin{cases}\label{density}
				c_{\gamma}(1-y^2)^{-\gamma},& y\in]0,1[;\\
				0,& \mbox{otherwise},
			\end{cases}
		\end{align}
		with $\gamma:=\frac{\alpha}{2(\alpha+2)}$ and $c_{\gamma}=\frac{\Gamma\left(\frac{3}{2}-\gamma\right)}{\sqrt{\pi}\Gamma(1-\gamma)}$;
		\item $Z$ and $Y$ are independent random variables;
		\item $\xi(t):=\frac{t^{\frac{\alpha}{2}+1}}{\frac{\alpha}{2}+1}$, $t\geq 0$.
	\end{itemize}
	This means that $u$ in \eqref{representation} corresponds to
	\begin{align}\label{integral rep}
		u(t,x)=\frac{t}{2\xi(t)}\int_{x-\xi(t)}^{x+\xi(t)}\varphi(y)c_{\gamma}\left(1-\frac{(x-y)^2}{\xi(t)^2}\right)^{-\gamma}dy,\quad t\geq 0,x\in\mathbb{R}.
	\end{align}
\end{theorem}	

We now investigate the behaviour of problem \eqref{tricomi} in the case of random singular initial condition, i.e.
\begin{align}
	\begin{cases}\label{tricomi stoch}
		U_{tt}(t,x)=t^{\alpha}U_{xx}(t,x),& t>0,x\in\mathbb{R};\\
		U(0,x)=0, U_t(0,x)=W(x),& x\in\mathbb{R},
	\end{cases}
\end{align}
where $\{W(x)\}_{x\in\mathbb{R}}$ stands for the so-called \emph{Gaussian white noise}. This means that $\{W(x)\}_{x\in\mathbb{R}}$ is defined to be the distributional derivative of a standard one dimensional two-sided Brownian Motion $\{B(x)\}_{x\in\mathbb{R}}$ defined on a reference probability space $(\Omega,\mathcal{F},\mathbb{P})$. We refer the reader to Section \ref{proof 2} for the definition and basic properties of the Gaussian white noise. 

The natural solution concept for investigating \eqref{tricomi stoch} is described as follows.  

\begin{definition}\label{def sol}
	We say that a random field $\{U(t,x)\}_{t\geq 0,x\in\mathbb{R}}$ solves problem \eqref{tricomi stoch} if for any sequence of mollifiers $\{\rho_n\}_{n\geq 1}$, i.e. for all $n\geq1$,
	\begin{itemize}
		\item $\rho_n\in C^{\infty}_c(\mathbb{R};[0,+\infty[)$ and $\int_{\mathbb{R}}\rho_n(y)dy=1$;
		\item the support of $\rho_n$ is contained in the closed ball $\overline{B(0,r_n)}$ for some non negative sequence $\{r_n\}_{n\geq 1}$ monotonically decreasing to zero,
	\end{itemize}
	the solution $\{U^{(n)}(t,x)\}_{t\geq 0,x\in\mathbb{R}}$ to 
	\begin{align}
		\begin{cases}\label{tricomi stoch molli}
			U^{(n)}_{tt}(t,x)=t^{\alpha}U^{(n)}_{xx}(t,x),& t>0,x\in\mathbb{R};\\
			U^{(n)}(0,x)=0, U^{(n)}_t(0,x)=W^{(n)}(x),& x\in\mathbb{R},
		\end{cases}
	\end{align}
	where
	\begin{align*}
		W^{(n)}(x):=\int_{\mathbb{R}}\rho_n(x-y)dB(y),\quad x\in\mathbb{R}
	\end{align*}
	converges in $\mathbb{L}^2(\Omega)$, as $n$ tends to infinity, to $\{U(t,x)\}_{t\geq 0,x\in\mathbb{R}}$.  
\end{definition}

We are now ready to state well posedness for the singular problem \eqref{tricomi stoch}; its proof can be found in Section \ref{proof 2}.

\begin{theorem}\label{main theorem}
	The Gaussian random field
	\begin{align}\label{sol tric stoch}
		U(t,x)=\frac{t}{2\xi(t)}\int_{x-\xi(t)}^{x+\xi(t)}c_{\gamma}\left(1-\frac{(x-y)^2}{\xi(t)^2}\right)^{-\gamma}dB(y),\quad t\geq 0,x\in\mathbb{R}.
	\end{align}
	with $\gamma$ and $\xi$ defined in Theorem \ref{th representation general}, is the unique solution in the sense of Definition \ref{def sol} to problem \eqref{tricomi stoch}.
\end{theorem}

\begin{remark}
	There are several results in the literature concerning the existence
	of solutions for the wave equation perturbed by additive noise in any
	space dimension.  We refer the reader to the papers \cite{Walsh} and
	\cite{Orsingher} for the case of one space dimension and to
	\cite{Dalang99}  and the contributions cited there for the
	multidimensional case. We remark that in \cite{Dalang99}  the author
	identifies a necessary and sufficient condition  relating additive noise terms to space dimension that allow for well posedness of the corresponding stochastic partial differential equations. \\
	The case of non strictly hyperbolic operators  on the contrary is not
	very much explored.  Here we mention some recent papers attempting to investigate how the degeneracies of the operator have an impact on stochastic perturbations: \cite{Ascanelli2,Ascanelli3,Ascanelli1,EnricoAlberto,EnricoLeonardo}. The results in this note are a new contribution in that direction.
\end{remark}

We now focus our attention on the problem
\begin{align}
	\begin{cases}\label{SPDE 2}
		V_{tt}(t,x)=t^{2}V_{xx}(t,x)-V_{x}(t,x),& t>0,x\in\mathbb{R};\\
		V(0,x)=0, V_t(0,x)=W(x),& x\in\mathbb{R}.
	\end{cases}
\end{align}
This corresponds, at the level of second order terms, to
\eqref{tricomi stoch} with $\alpha=2$ but differs from that for the
presence of a lower-order term, i.e. $V_x(t,x)$. Notice that, from the
point of view of the theory of hyperbolic operators with double
characteristics, the operators in \eqref{tricomi stoch} and
\eqref{SPDE 2} are the same: lower-order terms do not change the fact
that the Cauchy problem is $C^{\infty}$-well posed. They are the
prototype of the so called \emph{effectively hyperbolic} case: see e.g.
\cite{MR4718636} for a recent general review. However, equations
\eqref{tricomi stoch} and \eqref{SPDE 2} behave differently when the
initial data is Gaussian white noise. While \eqref{tricomi stoch}
admits a well defined square integrable solution, as stated in Theorem
\ref{main theorem}, this no longer happens for equation \eqref{SPDE 2}. For the proof of Theorem \ref{main theorem 2} see Section \ref{proof 2}. In the sequel we will implicitly adapt Definition \ref{def sol} with the straightforward modifications needed to treat the new operators under consideration.   

\begin{theorem}\label{main theorem 2}
	The Cauchy problem for the Tricomi-type partial differential equation with Gaussian white noise initial data \eqref{SPDE 2} doesn't admit any solution in the sense of Definition \ref{def sol}. 
\end{theorem}	

If we allow the initial data in \eqref{SPDE 2} to be slightly more regular than $\{W(x)\}_{x\in\mathbb{R}}$, more precisely a fractional Gaussian white noise with Hurst parameter $H\in (1/2,1)$, then the well posedness is restored. 	

\begin{corollary}\label{corollary}
	The problem   
	\begin{align}
		\begin{cases}\label{SPDE frac}
			V_{tt}(t,x)=t^{2}V_{xx}(t,x)-V_{x}(t,x),& t>0,x\in\mathbb{R};\\
			V(0,x)=0, V_t(0,x)=W_H(x),& x\in\mathbb{R}.
		\end{cases}
	\end{align}
	where $\{W_H(x)\}_{x\in\mathbb{R}}$ is now a fractional Gaussian white noise with Hurst parameter $H\in (1/2,1)$, admits a unique solution in the sense of Definition \ref{def sol}. 	
\end{corollary}

Proof of Corollary \ref{corollary} and definition of fractional Gaussian white noise with Hurst parameter $H\in (1/2,1)$ can be found in Section \ref{proof 2}.

We are now going to show that if the Tricomi-type operator is replaced by the classical wave operator then the introduction of a lower-order term doesn't modify the well posedness of the Cauchy problem with Gaussian white noise initial data. This confirms the existence of an interplay between degeneracy properties of differential operators and well posedness of some related stochastic Cauchy problems.  \\
Setting $\alpha=0$ in \eqref{tricomi stoch} one gets immediately from Theorem \ref{main theorem} that the unique random field solution of the problem 
\begin{align}
	\begin{cases}\label{wave}
		U_{tt}(t,x)=U_{xx}(t,x),& t>0,x\in\mathbb{R};\\
		U(0,x)=0, U_t(0,x)=W(x),& x\in\mathbb{R},
	\end{cases}
\end{align}
is given by
\begin{align}\label{solution wave}
	U(t,x)=\frac{1}{2}\int_{x-t}^{x+t}dB(y),\quad t\geq 0,x\in\mathbb{R};
\end{align}
since
\begin{align*}
	U(t,x)\in\mathbb{L}^2(\Omega)\quad\mbox{ for all }t\geq 0,x\in\mathbb{R},
\end{align*}
the solution \eqref{solution wave} is a well defined random field. We now perturb problem \eqref{wave} with a lower-order term, namely
\begin{align}
	\begin{cases}\label{wave +}
		V_{tt}(t,x)=V_{xx}(t,x)+V_x(t,x),& t>0,x\in\mathbb{R};\\
		V(0,x)=0, V_t(0,x)=W(x),& x\in\mathbb{R}.
	\end{cases}
\end{align}
A direct verification shows that its solution is given by
\begin{align*}
	V(t,x)=\frac{1}{2}\int_{x-t}^{x+t}e^{(y-x)/2}\mathtt{J}_0\left(\frac{1}{2}\sqrt{t^2-(x-y)^2}\right)dB(y),\quad t>0,x\in\mathbb{R};
\end{align*}
here $\mathtt{J}_0$ denotes the Bessel function of first kind of order zero. Computing the $\mathbb{L}^2(\Omega)$-norm of $V(t,x)$ we see that
\begin{align*}
	\mathbb{E}[V(t,x)^2]=\frac{1}{4}\int_{x-t}^{x+t}e^{y-x}\mathtt{J}_0\left(\frac{1}{2}\sqrt{t^2-(x-y)^2}\right)^2dy
\end{align*}
and the last integral turns out to be finite for all $t\geq 0$ and $x\in\mathbb{R}$ being the integral of a continuous function on a compact set. Therefore, both problems \eqref{wave} and \eqref{wave +} admit a unique well defined random field solution. 

The paper is organized as follows: in Section \ref{proof 1} we prove Theorem \ref{th representation general}. We take the form of the density of $Y$ in \eqref{representation} as given (being suggested by the results from \cite{Y}), for an exponent $\gamma$ determined by the proof, and search for a function $\xi$ that matches the probabilistic representation \eqref{representation} with the solution to the Cauchy problem \eqref{tricomi}. This approach leads to a system of non linear differential equations for $\xi$ where the existence of a solution is immediate but the uniqueness requires a fixed point argument whose details are presented in our proof. Then, Section \ref{proof 2} is devoted to the proofs of Theorem \ref{main theorem}, Theorem \ref{main theorem 2} and Corollary \ref{corollary}. This section also provides some basic definitions and properties of the Gaussian white noise and the fractional Gaussian white noise with Hurst parameter $H\in(1/2,1)$ which are taken as initial data of our Tricomi-type equations.

\section{Proofs of Theorem \ref{th representation general}}\label{proof 1}

Let $u$ be the function defined by
\begin{align}\label{b}
	u(t,x):=t\mathbb{E}[\varphi(x+\xi(t) Z Y)],\quad t\geq 0,x\in\mathbb{R},
\end{align}
with $Z$, $Y$ and $\gamma$ according to the statement of Theorem \ref{th representation general}. Our aim is to find a twice continuously differentiable function $\xi$ such that \eqref{b} solves the Cauchy problem \eqref{tricomi}.
First of all, since
\begin{align*}
	u_t(t,x)=\mathbb{E}[\varphi(x+\xi(t) Z Y)]+t\dot{\xi}(t)\mathbb{E}[\varphi'(x+\xi(t) Z Y)ZY],
\end{align*}
we see that the initial conditions in \eqref{tricomi} are satisfied if $\xi(0)=0$. Moreover,
\begin{align}\label{a}
	u_{tt}(t,x)=&\dot{\xi}(t)\mathbb{E}[\varphi'(x+\xi(t) Z Y)ZY]+(\dot{\xi}(t)+t\ddot{\xi}(t))\mathbb{E}[\varphi'(x+\xi(t) Z Y)ZY]\nonumber\\
	&+t\dot{\xi}(t)^2\mathbb{E}[\varphi''(x+\xi(t) Z Y)Z^2Y^2]\nonumber\\
	=&(2\dot{\xi}(t)+t\ddot{\xi}(t))\mathbb{E}[\varphi'(x+\xi(t) Z Y)ZY]+t\dot{\xi}(t)^2\mathbb{E}[\varphi''(x+\xi(t) Z Y)Y^2],	
\end{align}
where the last equality follows from $Z^2=1$. Let us now focus on the first expectation from the last member above:
\begin{align*}
	\mathbb{E}[\varphi'(x+\xi(t) Z Y)ZY]=&\mathbb{E}[\mathbb{E}[\varphi'(x+\xi(t) Z Y)ZY|\sigma(Z)]]\\
	=&\mathbb{E}[Z\mathbb{E}[\varphi'(x+\xi(t) Z Y)Y|\sigma(Z)]]\\
	=&\mathbb{E}\left[Z\int_0^1\varphi'(x+\xi(t) Z y)yc_{\gamma}(1-y^2)^{-\gamma}dy\right]\\
	=&-\frac{c_{\gamma}}{2}\mathbb{E}\left[Z\int_0^1\varphi'(x+\xi(t) Z y)\frac{d}{dy}\left(\frac{(1-y^2)^{-\gamma+1}}{-\gamma+1}\right)dy\right]\\	
	=&-\frac{c_{\gamma}}{2}\mathbb{E}\left[Z\left[\varphi'(x+\xi(t) Z y)\frac{(1-y^2)^{-\gamma+1}}{-\gamma+1}\right]_{y=0}^{y=1}\right]\\
	&+\frac{c_{\gamma}}{2}\xi(t)\mathbb{E}\left[Z^2\int_0^1\varphi''(x+\xi(t) Z y)\frac{(1-y^2)^{-\gamma+1}}{-\gamma+1}dy\right]\\
	=&\frac{c_{\gamma}}{2}\frac{\varphi'(x)}{-\gamma+1}\mathbb{E}\left[Z\right]\\
	&+\frac{\xi(t)}{2}\mathbb{E}\left[\int_0^1\varphi''(x+\xi(t) Z y)\frac{(1-y^2)}{-\gamma+1}c_{\gamma}(1-y^2)^{-\gamma+1}dy\right]\\
	=&\frac{\xi(t)}{2}\mathbb{E}\left[\varphi''(x+\xi(t) Z Y)\frac{(1-Y^2)}{-\gamma+1}\right].
\end{align*}
Here, in the second to last equality we utilized the fact that the function $y\mapsto (1-y^2)^{-\gamma+1}$ is equal to zero at $y=1$ (recall that $\gamma<1$) while in the last equality we exploited the identity $\mathbb{E}[Z]=0$. Therefore, we can continue the computation in \eqref{a} as follows:
\begin{align*}
	u_{tt}(t,x)=&(2\dot{\xi}(t)+t\ddot{\xi}(t))\mathbb{E}[\varphi'(x+\xi(t) Z Y)ZY]+t\dot{\xi}(t)^2\mathbb{E}[\varphi''(x+\xi(t) Z Y)Y^2]\\
	=&(2\dot{\xi}(t)+t\ddot{\xi}(t))\frac{\xi(t)}{2(-\gamma+1)}\mathbb{E}\left[\varphi''(x+\xi(t) Z Y)(1-Y^2)\right]+t\dot{\xi}(t)^2\mathbb{E}[\varphi''(x+\xi(t) Z Y)Y^2]\\
	=&(2\dot{\xi}(t)+t\ddot{\xi}(t))\frac{\xi(t)}{2(-\gamma+1)}\mathbb{E}\left[\varphi''(x+\xi(t) Z Y)\right]\\
	&+\left(t\dot{\xi}(t)^2-(2\dot{\xi}(t)+t\ddot{\xi}(t))\frac{\xi(t)}{2(-\gamma+1)}\right)\mathbb{E}[\varphi''(x+\xi(t) Z Y)Y^2].
\end{align*}
A comparison between the first and last members of the above chain of equalities gives
\begin{align*}
	u_{tt}(t,x)=&(2\dot{\xi}(t)+t\ddot{\xi}(t))\frac{\xi(t)}{2(-\gamma+1)}\mathbb{E}\left[\varphi''(x+\xi(t) Z Y)\right]\\
	&+\left(t\dot{\xi}(t)^2-(2\dot{\xi}(t)+t\ddot{\xi}(t))\frac{\xi(t)}{2(-\gamma+1)}\right)\mathbb{E}[\varphi''(x+\xi(t) Z Y)Y^2].
\end{align*}  
On the other hand, differentiating twice with respect to $x$ in \eqref{b} we obtain
\begin{align*}
	u_{xx}(t,x)=t\mathbb{E}[\varphi''(x+\xi(t) Z Y)].
\end{align*}
Therefore, if we want $u$ to solve
\begin{align*}
	u_{tt}(t,x)=t^{\alpha}u_{xx}(t,x)\quad\mbox{ for all }t>0,x\in\mathbb{R},
\end{align*}
we have to find a function $\xi$ such that
\begin{align}
	\begin{cases}\label{ODE}
		(2\dot{\xi}(t)+t\ddot{\xi}(t))\frac{\xi(t)}{2(-\gamma+1)}=t^{\alpha+1},& t>0;\\
		t\dot{\xi}(t)^2-(2\dot{\xi}(t)+t\ddot{\xi}(t))\frac{\xi(t)}{2(-\gamma+1)}=0,& t>0;\\
		\xi(0)=0
	\end{cases}
\end{align} 
It is clear, comparing the first to equations in \eqref{ODE}, that the function $\xi(t)=\frac{t^{\frac{\alpha}{2}+1}}{\frac{\alpha}{2}+1}$, $t\geq 0$ is a solution to \eqref{ODE}. On the other hand, due to the non linear character of the differential equation \eqref{ODE} uniqueness of solution is not straightforward. We now are going to prove that; first of all we rewrite the second equation in \eqref{ODE} like
\begin{align*}
	2(1 -\gamma)t \dot{\xi}^{2}=& 2\xi\dot{\xi} + t\xi\ddot{\xi}\\
	=&\frac{1}{2}\frac{d}{dt}\xi^{2} + \xi(\dot{\xi} + t\ddot{\xi} )\\
	=&\frac{1}{2}\frac{d}{dt}\xi^{2} + \xi\frac{d}{dt}(t\dot{\xi}),
\end{align*}
i.e.
\begin{align*}
	2(1 -\gamma)t \dot{\xi}^{2}=\frac{1}{2}\frac{d}{dt}\xi^{2} + \xi\frac{d}{dt}(t\dot{\xi}).
\end{align*}
Integrating both sides from $ 0 $ to $ t $ we get
\begin{align*}
	2(1 -\gamma)\int_{0}^{t}s\dot{\xi}^{2}(s)ds = \frac{1}{2}\xi^{2}(t)
	+ t\xi(t)\xi'(t) - \int_{0}^{t}s\dot{\xi}^{2}(s)ds,
\end{align*}	
that means
\begin{align*}
	(3 - 2 \gamma)\int_{0}^{t}s\dot{\xi}^{2}(s)ds =&
	\frac{1}{2}\xi^{2}(t) + \frac{t}{2}\frac{d}{dt}(\xi(t)^{2})\\
	=& \frac{1}{2}\frac{d}{dt}(t\xi(t)^{2}).
\end{align*}
We now integrate the last equation from $0$ to $T$, a fixed positive constant, and we have
\begin{align*}
	2(3 - 2 \gamma)\int_{0}^{T}\int_{0}^{t}s\dot{\xi}^{2}(s)ds dt = T\xi(T)^{2}
\end{align*}
that is
\begin{align}
	\label{eq:20}
	2(3 - 2 \gamma)\int_{0}^{T}(T-s)s\dot{\xi}^{2}(s)ds = T\xi(T)^{2},\quad \xi(0)=0.
\end{align}
Now we show that (\ref{eq:20}) has just one solution. We are interested in studying the non linear operator
\begin{align*}
	\mathcal{T}(g)(t) := \sqrt{ \frac{2L}{t} \int_0^t (t - s)s(g'(s))^2 ds}, \quad t\in [0,T],
\end{align*}
with $ L = 3 - 2\gamma $
and proving that the equation
\begin{align*}
	g(t) = \mathcal{T}(g)(t),\quad t\in [0,T]
\end{align*}
has a unique solution in a suitable function space. To this aim, let us define the Banach space
\begin{align*}
	X := \{ g \in C^1([0,T]) \mid g(0) = 0 \}, \quad \|g\| := \sup_{t \in [0,T]} |g'(t)|
\end{align*}
so that we control the derivative $g'$ (since the operator $\mathcal{T}$ depends on $g'$). It can be shown that $\mathcal{T}(g)(t) \in C^1([0,T])$ and that $\mathcal{T}(g)(0) = 0$, so that $\mathcal{T}$ maps $X$ into itself. Suppose $g_1, g_2 \in X$ and let $h(s) := g_1'(s) - g_2'(s)$. Define
\begin{align*}
	I_1 := \int_0^t (t-s)s (g_1'(s))^2 \, ds, \quad
	I_2 := \int_0^t (t-s)s (g_2'(s))^2 \, ds.
\end{align*}
Then,
\begin{align*}
	|\mathcal{T}(g_1)(t) - \mathcal{T}(g_2)(t)|
	&= \left| \sqrt{ \frac{2L}{t} I_1 } - \sqrt{ \frac{2L}{t} I_2 } \right| \\
	&= \sqrt{ \frac{2L}{t} } \cdot |\sqrt{I_1} - \sqrt{I_2}| \\
	&\leq \sqrt{ \frac{2L}{t} } \cdot \frac{|I_1 - I_2|}{\sqrt{I_1} + \sqrt{I_2}}.
\end{align*}
To estimate $|I_1 - I_2|$ we write
\begin{align*}
	|I_1 - I_2| 
	&= \left| \int_0^t (t-s)s \left[(g_1'(s))^2 - (g_2'(s))^2 \right] ds \right| \\
	&= \left| \int_0^t (t-s)s (g_1'(s) + g_2'(s))h(s) \, ds \right| \\
	&\leq 2M \|g_1 - g_2\| \int_0^t (t-s)s \, ds \\
	&= 2M \|g_1 - g_2\| \cdot \frac{t^3}{6},
\end{align*}
where $M := \max\{\|g_1\|, \|g_2\|\}$. Thus,
\begin{align*}
	|\mathcal{T}(g_1)(t) - \mathcal{T}(g_2)(t)| 
	\leq \sqrt{ \frac{2L}{t} } \cdot \frac{2M \|g_1 - g_2\| \cdot \frac{t^3}{6}}{\sqrt{I_1} + \sqrt{I_2}}.
\end{align*}
For small $t$, this implies that
\begin{align*}
	\|\mathcal{T}(g_1) - \mathcal{T}(g_2)\| \leq K T^{1/2} \|g_1 - g_2\|
\end{align*}
for some constant $K$ depending on $L$ and $M$. Therefore, for sufficiently small $T$, the operator $\mathcal{T}$ is a contraction and since
\begin{itemize}
	\item $(X, \|\cdot\|)$ is a Banach space;
	\item $\mathcal{T}$ is a contraction for small enough $T$,
\end{itemize}
we may apply the Banach fixed point theorem to conclude that $\mathcal{T}$ has a unique fixed point in $X$.\\
For larger values of $T$, one may apply the contraction mapping principle iteratively on subintervals (standard continuation method), or appeal to the general theory of non linear Volterra integral equations of the second kind. These are known to have unique solutions under continuity and Lipschitz-type conditions, both of which are satisfied here. The integro-differential equation
\begin{align*}
	t g(t)^2 = 2L \int_0^t (t - s)s (g'(s))^2ds,\quad t\in [0,T],
\end{align*}
has a unique solution $g(t) \in C^1([0,T])$ with $g(0) = 0$ for any
finite $T > 0$. \\
Since the integral representation \eqref{integral rep} follows immediately from \eqref{representation}, the proof of Theorem \ref{th representation general} is complete.

\section{Proofs of Theorem \ref{main theorem}, Theorem \ref{main theorem 2} and Corollary \ref{corollary}}\label{proof 2}

We start with a brief description of the Gaussian white noise $\{W(x)\}_{x\in\mathbb{R}}$. For more details we refer the reader to one of the books \cite{SPDEbook,Nualart}.\\
Let $\{W(\varphi)\}_{\varphi\in L^2(\mathbb{R})}$ be an isonormal Gaussian process, i.e. $\{W(\varphi)\}_{\varphi\in L^2(\mathbb{R})}$ is a family of centred Gaussian random variables defined on a common probability space $(\Omega,\mathcal{F},\mathbb{P})$ such that the map
\begin{align*}
	L^2(\mathbb{R})\ni\varphi\mapsto W(\varphi)
\end{align*}
is linear and
\begin{align*}
	\mathbb{E}[W(\varphi_1)W(\varphi_2)]=\int_\mathbb{R}\varphi_1(y)\varphi_2(y)dy.
\end{align*}	
In the sequel we will write
\begin{align*}
	B(x):=W(\mathtt{1}_{[0,x]}),\quad x\in\mathbb{R};
\end{align*}
this process possesses a continuous modification that turns out to be a one dimensional two-sided Brownian motion; with this notations we have
\begin{align*}
	\int_\mathbb{R}\varphi(y)dB(y)=W(\varphi),\quad \varphi\in L^2(\mathbb{R});
\end{align*}
this corresponds to the Wiener-It\^o integral of $\varphi$ with respect to the two-sided Brownian motion $\{B(x)\}_{x\in\mathbb{R}}$. The distributional derivative
\begin{align*}
	W(x):=\frac{\partial}{\partial x}B(x),\quad x\in\mathbb{R},
\end{align*}
is called \emph{Gaussian white noise} and will play the role of initial data in the partial differential equations under investigation. \\
Similarly, $\{W_H(\varphi)\}_{\varphi\in \mathcal{H}}$ is a fractional isonormal Gaussian process with Hurst parameter $H\in (1/2,1)$ if 
\begin{align*}
	\mathcal{H}:=\left\{\varphi:\mathbb{R}\to\mathbb{R}:\int_{\mathbb{R}^2}|\varphi(x)||\varphi(y)|H(2H-1)|x-y|^{2H-2}dxdy<+\infty\right\}
\end{align*}
and the family $\{W_H(\varphi)\}_{\varphi\in \mathcal{H}}$ is centred Gaussian and such that the map
\begin{align*}
	\mathcal{H}\ni\varphi\mapsto W_H(\varphi)
\end{align*}
is linear and
\begin{align*}
	\mathbb{E}[W_H(\varphi_1)W_H(\varphi_2)]=\int_{\mathbb{R}^2}\varphi_1(x)\varphi_2(y)|H(2H-1)|x-y|^{2H-2}dxdy.
\end{align*}	
In the sequel we will write
\begin{align*}
	B_H(x):=W_H(\mathtt{1}_{[0,x]}),\quad x\in\mathbb{R};
\end{align*}
this process possesses a continuous modification that turns out to be a one dimensional two-sided fractional Brownian Motion with Hurst parameter $H\in (1/2,1)$; with this notation we have
\begin{align*}
	\int_\mathbb{R}\varphi(y)dB_H(y)=W_H(\varphi),\quad \varphi\in \mathcal{H}.
\end{align*}
The distributional derivative
\begin{align*}
	W_H(x):=\frac{\partial}{\partial x}B_H(x),\quad x\in\mathbb{R},
\end{align*}
is called \emph{fractional Gaussian white noise}.

We now prove Theorem \ref{main theorem}. According to Definition \ref{def sol} to find the solution of \eqref{tricomi stoch} we have to investigate problem \eqref{tricomi stoch molli}. Thanks to Theorem \ref{th representation general} the solution to \eqref{tricomi stoch molli} is given by
\begin{align*}
	U^{(n)}(t,x)=\frac{t}{2\xi(t)}\int_{x-\xi(t)}^{x+\xi(t)}c_{\gamma}\left(1-\frac{(x-y)^2}{\xi(t)^2}\right)^{-\gamma}W^{(n)}(y)dy.
\end{align*}
On the other hand, recalling the definition of $W^{(n)}(y)$ and invoking Stochatic Fubini's theorem, see for instance \cite{DaPrato}, we can write
\begin{align*}
	U^{(n)}(t,x)&\quad=\frac{t}{2\xi(t)}\int_{x-\xi(t)}^{x+\xi(t)}c_{\gamma}\left(1-\frac{(x-y)^2}{\xi(t)^2}\right)^{-\gamma}\left(\int_{\mathbb{R}}\rho_n(y-z)dB(z)\right)dy\\
	&\quad=\int_{\mathbb{R}}\left(\frac{t}{2\xi(t)}\int_{x-\xi(t)}^{x+\xi(t)}c_{\gamma}\left(1-\frac{(x-y)^2}{\xi(t)^2}\right)^{-\gamma}\rho_n(y-z)dy\right)dB(z).
\end{align*}
Now, since
\begin{align*}
	&\int_{x-\xi(t)}^{x+\xi(t)}c_{\gamma}\left(1-\frac{(x-y)^2}{\xi(t)^2}\right)^{-\gamma}\rho_n(y-z)dy\\
	&\quad=\int_{\mathbb{R}}\mathtt{1}_{[x-\xi(t),x+\xi(t)]}(y)c_{\gamma}\left(1-\frac{(x-y)^2}{\xi(t)^2}\right)^{-\gamma}\rho_n(y-z)dy
\end{align*}
and the function
\begin{align*}
	y\mapsto \mathtt{1}_{[x-\xi(t),x+\xi(t)]}(y)c_{\gamma}\left(1-\frac{(x-y)^2}{\xi(t)^2}\right)^{-\gamma}
\end{align*}
belongs to $L^2(\mathbb{R})$ for all $t>0$ and $x\in\mathbb{R}$ (notice that $\gamma\in [0,1/2[$ for all $\alpha\geq 0$), by the approximation properties of mollifiers (see \cite{Evans}) we conclude that
\begin{align*}
	\int_{x-\xi(t)}^{x+\xi(t)}c_{\gamma}\left(1-\frac{(x-y)^2}{\xi(t)^2}\right)^{-\gamma}\rho_n(y-z)dy\to \mathtt{1}_{[x-\xi(t),x+\xi(t)]}(z)c_{\gamma}\left(1-\frac{(x-z)^2}{\xi(t)^2}\right)^{-\gamma}
\end{align*}
in $L^2(\mathbb{R})$ as $n$ tends to infinity and this in turn implies by virtue of the It\^o isometry $\mathbb{E}[W(\varphi)^2]=\int\limits_{\mathbb{R}}\varphi(y)^2dy$ that
\begin{align*}
	U^{(n)}(t,x)\to \frac{t}{2\xi(t)}\int_{x-\xi(t)}^{x+\xi(t)}c_{\gamma}\left(1-\frac{(x-y)^2}{\xi(t)^2}\right)^{-\gamma}(y)dB(y)
\end{align*}
in $\mathbb{L}^2(\Omega)$ as $n$ tends to infinity, thus proving Theorem \ref{main theorem}.

Next, we present the proof of Theorem \ref{main theorem 2}. We start by considering the Cauchy problem 
\begin{align}
	\begin{cases}\label{tricomi ++}
		v_{tt}(t,x)=t^{2}v_{xx}(t,x)-v_x(t,x),& t>0,x\in\mathbb{R};\\
		v(0,x)=0, v_t(0,x)=\varphi(x),& x\in\mathbb{R},
	\end{cases}
\end{align}
where $\varphi:\mathbb{R}\to\mathbb{R}$ is a twice continuously differentiable function. A direct verification shows that 
\begin{align}\label{v}
	v(t,x):=\int_0^t\varphi\left(x-\frac{t^2}{2}+z^2\right)dz,\quad t\geq 0,x\in\mathbb{R}
\end{align}
is the unique classical solution to \eqref{tricomi ++} and with a change of variable we rewrite \eqref{v} as
\begin{align}\label{v 2}
	v(t,x)=\int_{x-\frac{t^2}{2}}^{x+\frac{t^2}{2}}\varphi\left(y\right)\frac{1}{2\sqrt{y-x+\frac{t^2}{2}}}dy,\quad t\geq 0,x\in\mathbb{R}.
\end{align}
Now, according to Definition \ref{def sol} to find a solution to \eqref{SPDE 2} we have to investigate the problem
\begin{align}
	\begin{cases}\label{z}
		V^{(n)}_{tt}(t,x)=t^{2}V^{(n)}_{xx}(t,x)-V^{(n)}_{x}(t,x),& t>0,x\in\mathbb{R};\\
		V^{(n)}(0,x)=0, V^{(n)}_t(0,x)=W^{(n)}(x),& x\in\mathbb{R},
	\end{cases}
\end{align}
which by virtue of \eqref{v 2} has solution
\begin{align}\label{zz}
	V^{(n)}(t,x)&=\int_{x-\frac{t^2}{2}}^{x+\frac{t^2}{2}}\frac{1}{2\sqrt{y-x+\frac{t^2}{2}}}W^{(n)}(y)dy\nonumber\\
	&=\int_{\mathbb{R}}\left(\int_{x-\frac{t^2}{2}}^{x+\frac{t^2}{2}}\frac{1}{2\sqrt{y-x+\frac{t^2}{2}}}\rho_n(y-z)dy\right)dB(z),\quad t\geq 0,x\in\mathbb{R}.
\end{align}
Since the function
\begin{align}\label{za}
	y\mapsto \mathtt{1}_{[x-\frac{t^2}{2},x+\frac{t^2}{2}]}(y)\frac{1}{2\sqrt{y-x+\frac{t^2}{2}}}
\end{align}
doesn't belong to $L^2(\mathbb{R})$ (it does belong to $L^p(\mathbb{R})$ for any $p\in[1,2[$), we see that $V^{(n)}(t,x)$ doesn't converge in $\mathbb{L}^2(\Omega)$ and hence problem \eqref{SPDE 2} doesn't admit a solution in the sense of Definition \ref{def sol}. 

Lastly, we prove Corollary \ref{corollary}. We replace \eqref{z} with
\begin{align*}
	\begin{cases}
		V^{(n)}_{tt}(t,x)=t^{2}V^{(n)}_{xx}(t,x)-V^{(n)}_{x}(t,x),& t>0,x\in\mathbb{R};\\
		V^{(n)}(0,x)=0, V^{(n)}_t(0,x)=W_H^{(n)}(x),& x\in\mathbb{R},
	\end{cases}
\end{align*}
and formula \eqref{zz} with
\begin{align}\label{zq}
	V^{(n)}(t,x)&=\int_{x-\frac{t^2}{2}}^{x+\frac{t^2}{2}}\frac{1}{2\sqrt{y-x+\frac{t^2}{2}}}W_H^{(n)}(y)dy\nonumber\\
	&=\int_{\mathbb{R}}\left(\int_{x-\frac{t^2}{2}}^{x+\frac{t^2}{2}}\frac{1}{2\sqrt{y-x+\frac{t^2}{2}}}\rho_n(y-z)dy\right)dB_H(z),\quad t\geq 0,x\in\mathbb{R}.
\end{align}
Notice that the function \eqref{za} belongs to $\mathcal{H}$; in fact,
\begin{align*}
	\int_{x-\frac{t^2}{2}}^{x+\frac{t^2}{2}}\int_{x-\frac{t^2}{2}}^{x+\frac{t^2}{2}}\frac{1}{2\sqrt{y-x+\frac{t^2}{2}}}\frac{1}{2\sqrt{z-x+\frac{t^2}{2}}}|y-z|^{2H-2}dydz
\end{align*} 
is a Selberg-type integral \cite{Forrester} whose fineteness is guaranteed for $H\in (1/2,1)$. Moreover, since $L^{1/H}(\mathbb{R})\subset\mathcal{H}$ and 
\begin{align*}
	\int_{x-\frac{t^2}{2}}^{x+\frac{t^2}{2}}\frac{1}{2\sqrt{y-x+\frac{t^2}{2}}}\rho_n(y-z)dy\to \mathtt{1}_{[x-\frac{t^2}{2},x+\frac{t^2}{2}]}(y)\frac{1}{2\sqrt{y-x+\frac{t^2}{2}}}
\end{align*}
in $L^p(\mathbb{R})$ for all $p\in [1,2[$, such convergence holds also in $\mathcal{H}$. This in turn implies the convergence of \eqref{zq} to 
\begin{align*}
	\int_{x-\frac{t^2}{2}}^{x+\frac{t^2}{2}}\frac{1}{2\sqrt{y-x+\frac{t^2}{2}}}dB_H(y)
\end{align*}
in $\mathbb{L}^2(\Omega)$ for any $H\in(1/2,1)$. The proof is complete.

\bibliographystyle{plain}

\bibliography{tricomiARXIV}

\end{document}